\documentclass[a4paper,11pt]{article}
\usepackage{amsmath,amscd,amsfonts,amssymb,latexsym}

\newtheorem{theorem}[equation]{Theorem}

\newtheorem{corollary}[equation]{Corollary}
\newtheorem{remark}[equation]{Remark}

\newtheorem{formula}[equation]{Formula}

\def\Hom{{\rm Hom}}
\def\C{\mathbb{C}}
\def\Q{\mathbb{Q}}
\def\T{\mathbb{T}}

\def\P{\mathbb{P}}

\def\h{\Phi}
\def\CQ{{\widehat Q}}
\def\CX{{\widehat X}}
\def\CCQ{Q^*}
\def\CCX{X^*}

\def\qed{\hfill$\Box$\s}
\def\s{\vskip10pt}
\def\proof{\noindent{\it Proof.~~}}

\date{August 2014}

\begin{document}

\title{\bf Equivariant Hirzebruch class for quadratic cones via degenerations}
\author{Ma\l gorzata Mikosz\\
\small Warsaw University of Technology,\\
\small ul.~Koszykowa 75, 00-662, Warszawa, Poland\\
\small emmikosz@mini.pw.edu.pl \and
Andrzej Weber\thanks{Supported by NCN grant 2013/08/A/ST1/00804
}
\\
\small Department of Mathematics of Warsaw University\\
\small Banacha 2, 02-097 Warszawa, Poland\\
\small aweber@mimuw.edu.pl}

\maketitle
Let $X$ be a smooth algebraic variety and $Y$ a subvariety. The cohomology class of $Y$ in $H^*(X)$, i.e.~the Poincar\,e dual of the fundamental class of $Y$, does not change when we deform $Y$ in a flat manner. A more subtle cohomological invariant of $Y$ is the Hirzebruch class $$td_y(Y\to X)\in H_*(X)\otimes\Q[y]$$ defined in \cite{BSY}. A flat family member $Y_t$ can be thought of as a fiber of a function
$$X\times \C\supset W\stackrel{\pi}\longrightarrow \C\,.$$
The difference between the Hirzebruch class of the generic fiber and the Hirzebruch class of the special fiber is measured by the appropriate version of Milnor class, studied in \cite{CMSS} for hypersurfaces and in \cite{MSS} the general case. The same phenomenon happens for the equivariant Hirzebruch class developed in \cite{We3}, compare
also with \cite[Sec.4]{Oh} for the equivariant Hirzebruch class in the context
of quotient stacks. We fix our attention on the varieties with torus action. If we are interested in local invariants of singularities, we study the localization of the equivariant Hirzebruch class $td_y^\T(Y\to X)$ at a fixed point. The bottom degree of the Hirzebruch class is the equivariant fundamental class, also called the multi-degree of the variety. It does not change in the deformation class.
For example, let $\CQ_n\subset \C^n$ be the cone over a quadric in $\P^{n-1}$, in other words, $\CQ_n$ in some coordinates is described by the Morse function $\sum_{i=1}^n x_i^2$. Let $\T=\C^*$ act on $\C^n$ diagonally. Then $[\CQ_n]$ is equal to $2t$, with $$t=c_1(\C)\in H^*_\T(pt)\simeq H^*_\T( \C^n)\simeq \Q[t]\,,$$ the first Chern class of the standard weight one representation. Indeed $\CQ_n$ can be equivariantly degenerated to the sum of two transverse hyperplane. The difference of the Hirzebruch classes is supported by the singular locus of the special member of the family. In the  case of quadratic cones ($\CQ_n$ and intersection of planes) both varieties have only rational singularities,
therefore (\cite[Example 3.2]{BSY}) their Hirzebruch classes for $y=0$  are equal to the Todd classes constructed by Baum-Fulton-MacPherson. The Todd class of a hypersurface $H$ of an ambient manifold $M$ are expressed by the  class  $[H]$ and the Todd class of $M$, precisely $i_*td(H)=td(M)(1-e^{-[H]})$, where $i$ is the inclusion $i:H\hookrightarrow M$, see eg.~\cite[Th. 18.3(4)]{Fu}. One easily generalizes this formula in the equivariant setting.
  Hence the Todd classes of $\CQ_n$ and $\CX_n$ are equal. (Alternatively one can apply Verdier specialization argument, which implies that  the Todd class of singular
spaces is constant in flat families, \cite{Ve}.) It follows that full Milnor class is divisible by $y$.

We would like to present how  the equivariant Hirzebruch class degenerates for the cone singularities. Our work started when we tried to analyze the equivariant Hirzebruch class of the cone.
For the fixed dimension $n$ it is easy to compute the corresponding polynomial. From initial sequence of coefficients it was hard to guess a closed formula and,  for example, to prove a kind of positivity studied in \cite[\S13]{We3}. Applying the degeneration method we find an answer.
An interesting reciprocity happens. The difference between the Hirzebruch
classes of the projective quadric $Q_n$ and  two intersecting projective
hyperplanes $X_n$ is the Hirzebruch class of the complement of another
projective quadric multiplied by $y$:
\begin{equation}td^\T_y(Q_n)-td^\T_y(X_n) = y\cdot td^\T_y(\P^{n-3}\setminus Q_{n-2})\,\label{formulaglowna}\end{equation}
(Formula \ref{proj}).
In the non-equivariant context this result should follow for example from
\cite[Thm.1.4, Rem.1.5]{CMSS} and the methods of \cite[Sec.5]{PaPr}.
(as explained later in Remark \ref{recenzent}).
 In this paper we even prove more directly a corresponding
result for the equivariant Hirzebruch classes. Using induction we find the
equivariant Hirzebruch classes of $Q_n$ and $\CQ_n$.

Having in mind the expression for Chern-Schwartz-MacPherson class of smooth open varieties via logarithmic forms \cite{Al}, it is more natural to compute the Hirzebruch class of the complement $\CCQ_n=\C^n\setminus \CQ_n$.
For $n=2m$ we obtain the expression
$$(1 + y)^2T^2 \sum_{i=1}^m(-y)^{m - i}\frac{
    (1 +y T )^{2 i-2} }{(1 - T)^{2 i}}$$
 and  for $n=2m+1$
$$ (-y)^m  \frac{( y + 1) T }{1 - T}+(1 + y)^2T^2 \sum_{i=1}^m(-y)^{m - i}\frac{
    (1 + y T )^{2 i-1} }{(1 - T)^{2 i+1}}  \,. $$
Here $T=e^{-t}$ and the given expression is equal to the Hirzebruch class divided by the Euler class of $0\in \C^n$, that is  $eu(0)=t^n$. The formulas are understood as  elements of the completed $H^*_\T(\C^n)[y]$ and localized in $t$. This ring is isomorphic to the ring of Lautent series in $t$ and polynomials in $y$, i.e. $\Q[[t]][t^{-1},y]$. (We will omit the completion in our notation for cohomology.)
The formulas follow from Corollary \ref{expl} by the specialization $T_i$ to one.
Taking the limit $y\to -1$ with $T=e^{-(y+1)t}$ we obtain the expression for the Chern-Schwartz-MacPherson class of $\CCX_n$ in equivariant cohomology of $\C^n$:
for $n=2m$
$$\sum_{i=0}^{m-1} t^{2 i} (1 + t)^{2 (m -  i - 1)}$$
and for $n=2m+1$
$$t^{2 m}+\sum_{i=0}^{m-1}t^{2 i} (1 + t)^{2( m -  i - 1)} $$
which, as one can check, agrees with the invariant of a conical set introduced in \cite{AlMa}, compare \cite[\S8]{We1}.
We note that the quadratic cone appears as a  singularity of Schubert varieties: the quadric $Q_n$ can be considered as a homogenous space with respect to $SO(n)$ and the codimension one Schubert variety is isomorphic to the projective cone over $Q_{n-2}$. It would be interesting to examine singularities of Schubert varieties from the point of view of degenerations, having in mind the work on smoothability \cite{Co1,Co2} and intersection theory \cite{CoVa}.

The presented computation in fact is a baby example of what can happen.
The aim of the paper is to show a bunch of computation of the Hirzebruch class based on Localization Theorem \ref{locH}. The {\bf Formulas \ref{proj}, \ref{con}} and {\bf \ref{dope}} are the outcome. They show how Milnor class may be realized geometrically. We hope that these formulas will find generalizations for some class of degenerations  of Schubert varieties.\s\s

We would like to thank the Referee for very careful reading of the manuscript. He has suggested many important improvements. The meaningful Remark \ref{recenzent} is due to him.

\section{Hirzebruch classes of projective quadrics}

To understand systematically the situation we consider a bigger torus preserving the quadric. One has to distinguish between the cases of even and odd $n$.
Let us index the coordinates in $\C^{2m}$ by integer numbers from $-m$ to $m$ omitting $0$ and consider the quadratic form in $\C^{2m}$ given by the formula
$$\sum_{i=1}^mx_{-i}x_{i}\,.$$
For $\C^{2m+1}$ allow the index $0$ and fix the quadratic form
$$x_0^2+\sum_{i=1}^mx_{-i}x_{i}\,.$$
Let $Q_n\subset\P^{n-1}$ be the quadric defined by vanishing of the quadratic form. It is an invariant variety with respect to the torus $\T_m=(\C^*)^m$ action coming from the representation with weights (i.e. characters) $$(-t_m,-t_{-m+1},\dots,t_{m-1},t_m)$$ if $n=2m$ and $$(-t_m,-t_{-m+1},\dots,0,\dots,t_{m-1},t_m)$$  for $n=2m+1$. Consider the equivariant Hirzebruch class $$td_y^{\T_m}(Q_n\to \P^{n-1})\in H^*_{\T_m}(\P^{n-1})[y]$$ and compare it with the Hirzebruch class of degeneration $X_n$ of $Q_n$ given by the equation $x_{-m}x_m=0$. The variety $X_n$ is the sum of the  two coordinate planes. We think of $X_n$ as the special fiber for $\lambda=0$ of the equivariant family given by the equation $$\lambda \sum_{i=1}^{m-1}x_{-i}x_{i}+x_{-m}x_{m}\qquad \text{or}\qquad
\lambda \left(x_0^2+\sum_{i=1}^{m-1}x_{-i}x_{i}\right)+x_{-m}x_{m}\,.$$
We will show that the difference of the Hirzebruch classes is the Hirzebruch class of $\C^{n-2}\setminus Q_{n-2}$ multiplied by $y$, i.e.~Formula \ref{proj}.
\begin{remark}\label{recenzent} \rm Let us explain why the formula (\ref{formulaglowna}) holds in non-equivariant cohomology\footnote{This remark is due to the Referee}. In $H^*(\P^{n-1})$
$$td_y(Q_n)-td_y(X_n) = y \cdot td_y(\P^{n-3}\setminus Q_{n-2})$$
should follow from results and techniques a la
\cite[Thm.1.4, Rem.1.5]{CMSS} and \cite[Sec.5]{PaPr}:
$$g =\sum^{m-1}
_{i=1} x_{-i}x_i\quad \text{and} \quad f = x_{-m}x_m$$ are both sections of the line bundle
${\cal O}(2)$ on $\P^{n-1}$, with $Z' := \{g = 0\}$  and $Z := X_n = \{f = 0\}$ transversal in
a stratified sense. Let
$$p :{\mathcal Z} := \{\lambda g + f = 0\}\subset \P^{n-1} \times \C \to \C$$
be the projection onto the last variable $\lambda$. Then the vanishing cycles
$\phi_p(\Q_{\mathcal Z})$ are supported by the critical locus $\P^{n-3} = \{x_{-m} = 0 = x_m\} \subset
X_n = \{p = 0\}$ of $p$. Moreover, the restriction of these vanishing cycles
to $Z\cap Z' = Q_{n-3} \subset \P^{n-3}$ should be zero by the argument of \cite[Sec.5]{PaPr}
(or \cite[part a) of the proof of Prop.~4.1]{MSS}). Moreover, the corresponding
nearby cycles can be calculated in terms of the generic fiber $Q_n = \{p = 1\}$,
since $p$ is quasi-homogeneous (i.e. equivariant for a suitable $\C^*$-action).
Then the stated formula above follows from \cite[Thm.1.4, Rem.1.5]{CMSS},
with the factor $y$ equal to the (reduced) $\chi_y$-genus of the transversal Milnor
fiber of an $A_1$-singularity $z^2 + w^2 = 0 = x_{-m}x_m$ in $\C^2$.
This
remark would be an alternative proof of our formula provided that one developed the general theory of Milnor class in the equivariant case.\end{remark}

 It is more convenient to work with complements of the closed varieties from the beginning. We will give formulas for complements of the quadrics, since then the components have better geometric interpretation.
To make the notation easier we identify the equivariant cohomology with respect to $\T_m$ with the subspace of $H^*_{\T_{m+1}}(\P^{n-1})\simeq H^*_{\T_{m}}(\P^{n-1})\otimes \Q[t_{m+1}]$ given by $t_{m+1} = 0$ and omit the index $m$ in $\T_m$. Also we will omit the ambient space in the notation. This should not lead to a confusion; enlarging the ambient space results in introducing of the factor, which is the Euler class of the normal bundle.
 We will use this for example
for the inclusions
 $\iota:\P^{n-3}\to \P^{n-1}$  into the first coordinates and the corresponding inclusions of the affine spaces.
 For an isolated  fixed point $p\in Q_{n-2}\subset \P^{n-3}\subset\P^{n-1}$ the quotient
 $\frac{td_y^\T( Q_{n-2})_{|p}}{eu(p)}$ (where $eu(p)\in H^*_\T(pt)$ is the Euler class of the ambient tangent representation) does not depend on the ambient space. After these remarks about notation we state our first formula:
\begin{formula}\label{proj} Consider the complements of the quadrics $X_n'=\P^{n-1}\setminus X_n$ and $Q_n'=\P^{n-1}\setminus Q_n$. We have the equation  $$td_y^\T(X'_n)-td_y^\T(Q'_n)=y\,  td_y^\T( Q'_{n-2}) $$
in the equivariant cohomology $ H^*_\T(\P^{n-1})[y]$
for $n>2$. For the closed varie\-ties we have
$$td_y^\T(Q_n)-td_y^\T(X_n)=y\,  td_y^\T( Q'_{n-2})\,. $$
\end{formula}

\section{Topological and analytic localization theorems}

First let us note that equivariant cohomology is a homotopy invariant, for example for any $\T$-representation $V$ the restriction map $H^*_\T(V)\to H^*_\T(\{0\})$ is an isomorphism. Therefore we get for free
$H^*_\T(V)\stackrel{\simeq}{\to} H^*_\T(V^\T)$.
We need much stronger property of equivariant cohomology.
The main tool for computations is the Localization Theorem, see \cite[Ch.XII \S6]{Bo} or \cite{Qu}:
\begin{theorem}[Topological Localization Theorem] {\rm\cite[Theorem 4.4]{Qu}}\break
Assume either $X$ is a compact topological space or that $X$ is paracompact, $cd_\Q(X) < \infty$. Suppose a compact torus $\T$ acts on $X$ and the set of identity components of the isotropy groups of points of X is finite.
 Then the restriction map $$H^*_\T(X)\to H^*_\T(X^\T)$$
is an isomorphism after localization in the multiplicative system generated by nontrivial characters.\label{locH} \end{theorem}
We apply Topological Localization Theorem to algebraic varieties with algebraic torus action. The fixed points of the compact torus are the same as the fixed points of the full torus. The theorem may be applied to {\it any} algebraic variety, but it may very well happen (exactly when $X^\T=\emptyset$) that the localized equivariant cohomology is trivial.

For differential manifolds the isomorphism was made explicit by Atiyah-Bott  and Berline-Vergne, see also \cite{EdGr}.
\begin{theorem}[Topological Localization Theorem] {\rm   \cite[page 9]{AB}, \cite{BV}} Let $\T$ be a compact torus and let $M$ be a compact $\T$-ma\-nifold. Let $$M^T=\bigsqcup_{\alpha\in I} F_\alpha$$ be the decomposition of the fixed point set into connected components. Denote by $\iota_\alpha:F_\alpha\to M$ the inclusion. Let $$eu(F_\alpha)\in H^*_\T(F_\alpha)\simeq H^*(F_\alpha)\otimes H^*_\T(pt)$$ be the equivariant Euler class of the normal bundle to $F_\alpha$. Let $S$ be the multiplicative system generated by nontrivial characters. Then\begin{enumerate}
 \item The class $eu(F_\alpha)$ is invertible in $S^{-1}H^*_\T(F_\alpha)$.
\item For any equivariant cohomology class $\omega\in H^*_\T(M)$ the equality holds
\begin{equation}\omega=\sum_{\alpha\in I}\iota_{\alpha *}\left(\frac{\iota_\alpha^*(\omega)}{eu(F_\alpha)}\right)\label{skladniki}\end{equation}
in $S^{-1}H^*_\T(M)$.
\end{enumerate}\end{theorem}
The resulting integration formula follows, \cite[Formula 3.8]{AB}.

The case of compact algebraic smooth varieties is special. The equiva\-riant cohomology with respect to an algebraic torus action is always a free module over $H^*_\T(pt)$ (see \cite{GKM} and the references therein). Therefore the restriction map
\begin{equation*}H^*_\T(M)\to H^*_\T(M^T)\end{equation*}
is a monomorphism. {\it The equality of the classes restricted to the fixed point set implies their equality.} We will use just this principle. Nevertheless, having in mind the formula (\ref{skladniki}), it is natural and convenient to consider the {\it localized Hirzebruch class} $$\frac{\iota_\alpha^*(td^\T_y(-))}{eu(F_\alpha)}$$ in the localized  cohomology of fixed point set components. The spaces we consider here have only isolated fixed point sets, thus the localized Hirzebruch classes are polynomials in $y$ with coefficients in the ring of Laurent polynomials in $t_i$'s. In fact the coefficients are rational functions in $T_i=e^{-t_i}$.
\section{Properties of equivariant Hirzebruch class}
Now we would like to recall basic properties of the equivariant Hirzebruch class, which in fact formally do not differ from the properties of the non-equivariant class. For an equivariant line bundle $L$ the class $td^\T_y(L)$ is given in equivariant cohomology by the power series $$t\frac{1 + y\,e^{-t}}{1- e^{-t}}\,,$$ with
$t$ the first equivariant Chern class of $L$. Then the corresponding
class of a vector bundle is given in terms of Chern roots, and the
class for a smooth manifold $M$ is the corresponding class of the tangent bundle $TM$.
In the localized classes of a smooth manifold appears then the
(corrected) factor $$\h(T) = \frac{1+y T}{1-T}$$
with $T = e^{-t}$ at the normal directions to the fixed point set.

The important properties of the equivariant Hirzebruch classes of singular varieties used
in this paper are:\begin{enumerate}
\item
the normalization for smooth spaces (the Hirzebruch class is a series in equivariant Chern classes of tangent bundle),
\item covariant functoriality under proper maps,
\item additivity.
\end{enumerate}
For example: Let $\pi:\widetilde M\to M$ be an equivariant proper morphism, with $\pi_{|\widetilde M\setminus E}$ an isomorphism on the image for some $E \subset M$, a closed invariant subspace (for example the blowup of the origin in $M = \C^n$ with
$E=\P^{n-1}$ the exceptional divisor, as used later on). Then
$$\pi_*(td_y^\T(\widetilde M)-td_y^\T(E)) = td_y^\T(M) - td_y^\T(\pi(E))\,. $$
As an example for additivity (or the inclusion-exclusion principle)
one can calculate:
$$td_y^\T(X_n) = td_y^\T(\{x_{m} = 0\})+td_y^\T(\{x_{-m} = 0\})-td_y^\T(\{x_{-m} =x_{m} = 0\})\,;$$
since $X_n= \{x_m = 0\}\cup   \{x_{-m} = 0\}$ but the intersection is counted twice.  In particular  one can calculate
in this simple way the class of the singular space $X_n$ in terms of classes
of smooth spaces.
\s
The next property follows from 1.-3.:
\begin{description}\item {4.} multiplicativity and, more generally, contravariant functoriality with respect to fibrations.\end{description}
For example if $p:\nu\to X$ is an equivariant vector bundle, then the Hirzebruch class of the total space of $\nu$ is equal to
\begin{equation}td_y^\T(Tot(\nu))=p^*\left(td_y^\T(\nu)\cdot td_y^\T(X)\right)\,.\label{covariant}\end{equation}
Here $td_y^\T(\nu)$ is understood as a characteristic class of a vector bundle.

\section{Proof of Formula \ref{proj}.} By Localization Theorem \ref{locH} it is enough to check equality at each fixed point of $\T$-action.
The fixed points $p_i$ corresponds to the coordinate lines in $\C^{n}$. Let us show the calculation for  even $n=2m$. At the point $p_
i$ the quadric is given by the equation $$u_{-i}+\sum_{j\not=i}u_{-j}u_{j}=0$$
in coordinates $u_j=x_j/x_i$. For a fixed point $p_i$ the Hirzebruch class $td_y^\T(Q_n)$ divided by Euler class of at ${p_i}$  (i.e.~the localized Hirzebruch class) is equal to the product
$$\frac1{eu(p_i)}td_y^\T(Q_n)=\prod_{\text{weights of }T_{p_i}Q_n} \h(e^{-w})\,.$$
Here the product is taken with respect to
the weights appearing in the tangent representation $T_{p_i}Q_n$ (see \cite[\S1]{We3}).  Let us set $t_{-i}=-t_i$ and $T_i=e^{-t_i}$. The weights of the tangent representation $T_{p_i}\P^{n-1}$ are equal to $t_j-t_i$ for $j\not=i$. The normal direction has weight $t_{-i}-t_i=-2t_i$.
   Since $Q'_{2m}$ is the complement of $Q_{2m}$ in $\P^{n-1}$, one gets by additivity
that
$$\frac 1{eu(p_i)}td_y^{\T}(Q'_{2m})_{p_i}=\frac 1{eu(p_i)}td_y^{\T}(\P^{n-1})_{p_i}-\frac 1{eu(p_i)}td_y^{\T}(Q_{2m})_{p_i}.$$
\begin{itemize}\item At each point $p_i$, $|i|\leq m$ the localized Hirzebruch class is equal to
$$ (\h(T_i^{-2})-1)\cdot\prod_{j=1,j\ne i}^{m} \h(T_jT_i^{-1}) \h(T_j^{-1}T_i^{-1})\,.$$
\end{itemize}
The class $td_y^{\T}(X'_{2m})$ is equal to
$$td_y^\T(\P^{n-1})-td_y^\T(\{x_m=0\})-td_y^\T(\{x_{-m}=0\})+td_y^\T(\{x_m=x_{-m}=0\})
\,.$$
Therefore the localized class $\frac1{eu(p_i)}td_y^{\T}(X'_{2m})_{|p_i}$
 is the following

  \begin{itemize}
\item at the points $p_i$, $|i|<m$
$$ \h(T_i^{-2})\cdot (\h(T_mT_i^{-1})-1)(\h(T_m^{-1} T_i^{-1})-1)\cdot\prod_{j=1,j\ne i}^{m-1} \h(T_jT_i^{-1}) \h(T_j^{-1}T_i^{-1})$$
since
$$\h(T_mT^{-1}_i)\h(T_m^{-1}T^{-1}_i)-\h(T_mT^{-1}_i)-\h(T_m^{-1}T^{-1}_i)+1=$$
$$=(\h(T_mT^{-1}_i)-1)(\h(T_m^{-1}T^{-1}_i)-1)$$
\item at the point $p_{i}$, $|i|=m$
$$ (\h(T_i^{-2})-1)\cdot\prod_{j=1}^{m-1} \h(T_jT_i^{-1}) \h(T_j^{-1}T_i^{-1})\,.$$

\end{itemize}
For the points $p_{-m}$ and $p_{m}$ which do not belong to $\iota(Q_{n-2})$ the considered classes are equal. At the   point $p_{i}$ for $|i|<m$ the classes $td_y^{\T}(Q'_{2m})$, $y\,td_y^{\T}(Q'_{2m-2})$ and $td_y^{\T}(X'_{2m})$ have the common factor
 $$\prod_{j=1,j\ne i}^{m-1} \h(T_jT_i^{-1}) \h(T_j^{-1}T_i^{-1})$$
  and
   it is enough to check the equality
\begin{align*}
\h(T_i^{- 2})\cdot (\h(T_mT_i^{-1} ) - 1)&\cdot (\h(T_m^{-1}T_i^{-1} ) - 1) -\\
-(\h(T_i^{-2}) -& 1) \cdot \h(T_mT_i^{-1} )\cdot \h( T_m^{-1}T_i^{-1} )\, =\,y (\h(T_i^{-2}) - 1)\,.\end{align*}
\s
\noindent After multiplying by $$(1-T_i^{- 2})\cdot(1-T_i^{-1} T_m)\cdot (1-T_i^{-1} T_m^{-1})$$ the equality reduces to
\begin{align*}
(1+yT_i^{- 2})\cdot (y+1)(T_mT_i^{-1} )\cdot  (y+1)( T_m^{-1}&T_i^{-1})  -\\
-(y+1)(T_i^{-2})\cdot (1+yT_m&T_i^{-1} )\cdot (1+y T_m^{-1}T_i^{-1} )=\\
=\,y (y+1)&(T_i^{-2})\cdot(1-T_mT_i^{-1} )\cdot (1- T_m^{-1}T_i^{-1}) \,,\end{align*}

\noindent which one verifies easily.
The proof for $n$ odd is identical except that all the expressions are multiplied by $\h(T_i^{\pm1})$.
\qed

Also for $n=2$ if we admit that $Q_0=\P^{-1}=\emptyset$ and $td_y(\emptyset)=0$ the Formula \ref{proj} holds.

\section{Affine cones}
Let us extend the torus action by adding one factor to $\T$. Now we consider $\T=(\C^*)^{m+1}$ the character of the additional coordinate of $\T$ is denoted by $t$ and $T=e^{-t}$. The weights of the action on $\C^n$ are $$(t+t_{-m},t-t_{-m+1},\dots,t+t_{m-1},t-t_m)$$ in the even case and $$(t+t_{-m},t-t_{-m+1},\dots,t,\dots,t+t_{m-1},t-t_m)$$ in the odd case. It does not change the action on $\P^{n-1}$ on which the additional  coordinate of $\T$ acts trivially.

\begin{formula}\label{con} Consider the complements of the affine cones
$\CCQ_n=\C^{n}\setminus \CQ_n$ and $\CCX_n=\C^{n}\setminus \CX_n$. In the equivariant cohomology $H^*_\T(\C^n)[y]$ we have the equation
 $$td_y^\T( \CCX_n)-td_y^\T(\CCQ_n)=y\, td_y^\T( \CCQ_{n-2})$$
 for $n\geq2$.
\end{formula}

\proof Let $Y$ denote $X_n$, $Q_n$ or $Q_{n-2}$.
Let $\pi:\widetilde \C^n\to \C^n$ be the blowup at the origin with  $i: \P^{n-1}\hookrightarrow \widetilde\C^n$
the
inclusion of the exceptional divisor. The  Hirzebruch class of $Y^*\subset \C^n$ can be computed by push-forward of the class
$td_y^\T(\pi^{-1}( Y^*))$, since $\pi:\widetilde \C^n\setminus \P^{n-1}  \to \C^n\setminus\{0\}$ is an isomorphism. Here we are using  functoriality of the equivariant Hirzebruch classes. The projection $p:\widetilde \C^n\to \P^{n-1}$ has a structure of a vector bundle $\nu={\cal O}(-1)$. We apply the formula (\ref{covariant}) and additivity for $\pi^{-1}(Y^*)=p^{-1}(Y')\setminus i(Y')$:
$$td_y^\T(Y^*)=\pi_*td_y^\T(\pi^{-1}( Y^*))=\pi_*p^*\left((td_y^\T(\nu)-c_1(\nu))\cdot td_y^\T( Y')\right)\,.$$
The expression is linear with respect to $td_y^\T( Y')$.
It follows that the li\-near relation (Formula \ref{proj}) among Hirzebruch classes $td_y^\T( X'_n)$, $td_y^\T(Q'_n)$ and $td_y^\T( Q'_{n-2})$ in $H^*_\T(\P^{n-1})[y]$ implies the corresponding relation in $H^*_\T(\C^{n})[y]$.
\qed

\begin{remark}\rm  More generally for the degeneration
$$\lambda \sum_{i=1}^{k}x_{-i}x_{i}+\sum_{i=k+1}^m x_{-i}x_{i}$$
(and similarly for $n$ odd) we have
$$td_y^\T(\CCQ_n)-td_y^\T(Y^*)=(-y)^{m-k}\, td_y^\T(\CCQ_{2k})$$ where $Y$ is the hypersurface  corresponding to $\lambda=0$. The general case follows from the case $k=m-1$, which was studied here.
\end{remark}

We obtain the explicit formula
\begin{corollary}\label{expl}
 The equivariant Hirzebruch class of the complement of the quadratic cone $\CCQ_n=\C^n\setminus\CQ_n$ is equal to\hfill\break
for $n=2m$
$$td_y^{\T}(\CCQ_n )=\sum_{k=0}^{m-1}(-y)^k
td_y^{\T}(\CCX_{n-2k} )$$
for $n=2m+1$
$$td_y^{\T}(\CCQ_n )=\sum_{k=0}^{m-1}(-y)^k
td_y^{\T}(\CCX_{n-2k} )+(-y)^m td_y^{\T}(\C\setminus 0)\,,$$
where
$$\frac{td_y^{\T}(\CCX_{2m} )}{eu(0)}=(\h(TT_{m})-1)\cdot(\h(TT_{m}^{-1})-1)
\cdot\prod_{j=0}^{m-1}\h(TT_j)\h(TT_j^{-1})$$  and
$$\frac{td_y^{\T}(\CCX_{2m+1})}{eu(0)}=\h(T)\cdot(\h(TT_{m})-1)\cdot(\h(TT_{m}^{-1})-1)
\cdot\prod_{j=0}^{m-1}\h(TT_j)\h(TT_j^{-1})\,,$$
$$\frac{td_y^{\T}(\C\setminus 0)}{eu(0)}=\h(T)-1\,.$$
\end{corollary}\s

\section{Positivity}
Now we will show that the Hirzebruch classes of $\CQ_n$ and $\CCQ_n$ satisfy certain positivity condition.
For a weight $w\in \Hom(\T,\C^*)$ let us set a new variable $S_w=e^{-w}-1$. Also let us set $\delta=-1-y$.

\begin{corollary}\label{pos1}
The Hirzebruch class of the complement of the affine cone $\CCQ_n=\C^n\setminus \CQ_n$ is equal to  a  polynomial in $\delta$ and  $S_w$  with nonnegative coefficients divided by the product of the variables $S_w$, where $w$ are the weights of the representation $\C^n$.\end{corollary}

\proof It suffices to note that for the standard  action of one dimensional torus on $\C$ we have (with $T=e^{-t}$ as before)
$$\frac{td_y^\T(\C)}{eu(0)}= \h(T)=\frac{1-T+(1+y)(T-1+1)}{1-T}=\frac{
S_t +\delta (S_t+1)}{S_t}$$
and
$$\frac{td_y^\T(\C\setminus \{0\})}{eu(0)}=\frac{td_y^\T(\C)}{eu(0)}-\frac{td_y^\T( \{0\})}{eu(0)}=\h(T)-1  =\frac{\delta  (S_t+1)}{S_t}.$$
Moreover, since $\widehat X'_{n-2k}=\C^{n-2-2k}\times (\C^*)^2$ for $k=0,\dots,m-1$, by multiplicativity,  the Hirzebruch class $td_y^\T(\widehat X_n')$ is a nonnegative expression. The claim for $\CCQ_n$ follows from
 Corollary \ref{expl}.
\qed

For the original closed varieties we have:

\begin{formula}\label{dope} $$td_y^{\T}(\CQ_n)-td_y^\T( \CX_n)=y\left( td_y^\T(\C^{n-2})- td_y^\T(\CQ_{n-2})\right)\,.$$
\end{formula}

\proof We rewrite the Formula  \ref{con} passing to the complement
$$\left(td_y^\T(\C^n)- td_y^\T(\CX_n)\right)-\left(td_y^\T(\C^n)- td_y^\T(\CQ_n)\right)
=y\, \left(td_y^\T(\C^{n-2})-td_y^\T( \CQ_{n-2})\right)\,.$$
Hence we obtain what is claimed.\qed

\begin{corollary}\label{pos2}
The Hirzebruch class of the affine cone of $\CQ_n$ is equal to  a  polynomial in $\delta$ and  $S_w$  with nonnegative coefficients divided by the product of the variables $S_w$, where $w$ are the weights of the representation $\C^n$.\end{corollary}
\proof
Transforming the Formula \ref{dope} we obtain that
$$td_y^{\T}(\CQ_n)=-y\,td_y^\T( \CQ_{n-2})+\left(td_y^\T( \CX_n)+y\,td_y^\T(\C^{n-2})\right)$$
\begin{equation}=-y\,td_y^\T( \CQ_{n-2})+td_y^{\T}(\C^{n-2})\cdot\frac{  -(1 + y)(T^2-1)  }{(1 - TT_m^{-1}) (1 - T T_m)}\,.\label{posdomcz}\end{equation}
Here we use additivity and multiplicativity of the
Hirzebruch class applied to the decomposition $\widehat X^n=\C^{n-2}\times (\C_+\cup\C_-\setminus \{0\})$ with $$\frac{td^\T_y(\C_\pm)}{eu(0)}=\frac{1+y TT^{\pm1}_ m}{1-TT^{\pm1} _m}\,.$$
The formula (\ref{posdomcz}) follows from the identity
$$\frac{1+y TT_ m}{1-TT _m}+\frac{1+y TT^{-1}_ m}{1-TT^{-1} _m}-1+y=\frac{  -(1 + y)(T^2-1)  }{(1 - TT_m^{-1}) (1 - T T_m)}\,.$$
We note that
$$\frac{- (1 + y)( T^{2}-1)}{(1-TT_m) (1- T T_m^{-1})}=\frac{\delta(S_t^2+2 S_t)}{S_{t+t_m}S_{t-t_m}}$$
is a positive expression. We proceed inductively having in mind that the coefficient before $td_y^\T( \CQ_{n-2})$ is $-y=1+\delta$.\qed

The Corollaries \ref{pos1} and \ref{pos2} confirm the general rule (not proved so far) that the local Hirzebruch classes of Schubert cells are positive expressions in the variables associated with tangent weights.

\end{document}